\documentclass[11pt,oneside]{article}
\usepackage{color}
\usepackage{psfrag}
\usepackage{epsfig}
\usepackage{graphicx}
\usepackage{amsmath}
\usepackage{amssymb}
\usepackage{amsthm}
\usepackage{cite}
\usepackage[latin1]{inputenc}
\usepackage{cancel}

\newtheorem{theorem}{Theorem}[section]

\newtheorem{remark}[theorem]{Remark}
\begin{document}

\title{On permutation trinomials of type $x^{2p^s+r}+x^{p^{s}+r} +\lambda x^r$}
\date{}
\author{Daniele Bartoli, Giovanni Zini
\thanks{The research of D. Bartoli and G. Zini
was supported in part by Ministry for Education, University
and Research of Italy (MIUR) (Project PRIN 2012 ``Geometrie di Galois e
strutture di incidenza'' - Prot. N. 2012XZE22K$_-$005)
 and by the Italian National Group for Algebraic and Geometric Structures
and their Applications (GNSAGA - INdAM).}
}

\maketitle

\begin{abstract}

We determine all permutation trinomials of type $x^{2p^s+r}+x^{p^s+r}+\lambda x^r$ over the finite field $\mathbb{F}_{p^t}$ when $(2p^s+r)^4<p^t$. This partially extends a previous result by Bhattacharya and Sarkar in the case $p=2$, $r=1$.
\end{abstract}

{\bf Keywords:} Permutation trinomials, exceptional polynomials.

\section{Introduction}

Let $p$ be a prime number, $t$ be a positive integer, $\mathbb{F}_{p^t}$ be the finite field with $p^t$ elements, and $f(x)$ be a polynomial over $\mathbb F_{p^t}$.
If $f:x\mapsto f(x)$ is a permutation of $\mathbb F_{p^t}$, then $f(x)$ is a \emph{permutation polynomial} (PP) of $\mathbb F_{p^t}$.
If $f(x)$ is a PP of $\mathbb F_{(p^t)^m}$ for infinitely many $m$, then $f(x)$ is an \emph{exceptional polynomial} over $\mathbb F_{p^t}$.
If both $f(x)$ and $f(x)+x$ are PPs of $\mathbb F_{p^t}$, then $f(x)$ is a complete permutation polynomial (CPP) of $\mathbb F_{p^t}$.

The study of permutation polynomials over finite fields is motivated not only by their theoretical importance, but also by their remarkable applications to cryptography, combinatorial designs, and coding theory; see for instance \cite{LC2007,MN1987,SGCGM2012}.
For a detailed introduction to old and new developments on permutation polynomials, see the survey \cite{H2015} and the references therein.
Permutation polynomials of monomial and binomial type have been intensively investigated, while much less is known on permutation trinomials; see \cite{DLWYY2015,H2015tri}.

In this note we characterize a certain class of permutation trinomials.
Let $s$ and $r$ be non-negative integers. For $\lambda\in\mathbb F_{p^t}$, denote by $f_\lambda(x)$ the polynomial
$$ f_\lambda(x)=x^{2p^s+r}+x^{p^s+r}+\lambda x^r \in \mathbb{F}_{p^t}[x]. $$
If $r=0$, define $d=0$. If $r\ne0$, write $r=p^u v$ with $u\geq0$ and $p\nmid v$, and define $d=2p^{s-u}+v$ if $u\leq s$, $d=2+p^{u-s}v$ if $u>s$; that is,
$$ d=(2p^s+r)/p^m,\quad m:=\max\{n\geq0 \;:\; p^n\mid(2p^s+r),\;p^n\mid(p^s+r),\;p^n\mid r \} . $$
We prove the following result.
\begin{theorem}\label{result}
Assume that $d^4<p^t$.
Then $f_{\lambda}(x)$ is a PP of $\mathbb F_{p^t}$ if and only if one of the following cases holds:
\begin{itemize}
\item $p=2$, $t$ is odd, and $f_{\lambda}(x)=x^3+x^2+x$ or $f_{\lambda}(x)=x^5+x^3+x$;
\item $p\equiv2\pmod3$, $t$ is odd, and $f_{\lambda}(x)=x^3+x^2+\frac{1}{3}x$.
\end{itemize}
\end{theorem}

The case $p=2$ and $r=1$ was already considered by Bhattacharya and Sarkar \cite{BS2017}, where the result proved for $f_\lambda(x)$ was then used to characterize permutation binomials of $\mathbb F_{p^{2tp^s}}$ of type $g_b(x)=x^{\frac{p^{2tp^s}-1}{p^t-1}+1}+bx$.
Here for $p>2$ and $r=1$ we go the opposite direction, using the characterization in \cite{BGQZ2017} for permutation binomials of type $g_b(x)$ to deduce the result for $f_\lambda(x)$.

Every permutation polynomial of $\mathbb F_{p^t}$ with degree less than $(p^t)^{1/4}$ is exceptional over $\mathbb F_{p^t}$; thus, the condition $d^4<p^t$ allows us to consider only exceptional polynomials.
For $r>1$, this leads to the non-existence of permutation trinomials of type $f_\lambda(x)$.

\section{Proof of Theorem \ref{result}}

Since the maps $x\mapsto x^{p^u}$ and $x\mapsto x^{p^s}$ are permutations of $\mathbb{F}_{p^t}$, we can assume that $u=0$ if $u\leq s$, and $s=0$ if $u>s$.
\begin{itemize}
\item Case $r=0$. Since $f_\lambda(x)=f_\lambda(-x-1)$, $f_\lambda(x)$ is not a PP of $\mathbb{F}_{p^t}$.
\item Case $r=1$ and $p=2$. The claim is proved in \cite[Theorem 1.3]{BS2017}.
\item Case $r=1$ and $p>2$.

Assume first that $s=0$, so that $f_{\lambda}(x)=x^3+x^2+\lambda x$.
By direct computation,
$$ \frac{f_\lambda(x)-f_\lambda(y)}{x-y}=x^2+xy+y^2+x+y+\lambda$$
splits into two linear components if and only if $p\neq 3$ and $\lambda=\frac{1}{3}$. In this case 
$$ \frac{f_\lambda(x)-f_\lambda(y)}{x-y}=\left(x+\frac{1+\sqrt{-3}}{2}y+\frac{1}{2}+\frac{\sqrt{-3}}{6}\right)\cdot\left(x+\frac{1-\sqrt{-3}}{2}y+\frac{1}{2}-\frac{\sqrt{-3}}{6}\right)$$
and the two components are not defined over $\mathbb{F}_{p^t}$ if and only if $-3$ is a non-square in $\mathbb F_{p^t}$. From \cite[Lemma 4.5]{G2002}, this is equivalent to require $t$ odd and $p\equiv2\pmod3$.

Now assume that $s>0$.
Let $\mu\in\mathbb F_{p^t}$ with $\mu^{p^s}=\lambda$, so that $f(x)=x\left(x^2+x+\mu\right)^{p^s}$. Let $b,b^{p^t}\in\mathbb F_{p^{2t}}$ be the zeros of $x^2+x+\mu$; then, for any $x\in\mathbb F_{p^t}$,
$$ f_{\lambda}(x)=x\left(x+b\right)^{p^s}(x+b^{p^t})^{p^s}=x\left(x+b\right)^{(p^t)^{2p^s-1}+(p^t)^{2p^s-2}+\cdots+p^t+1} .$$
Suppose that $f_\lambda(x)$ is a PP of $\mathbb F_{p^t}$; in particular, $b\notin\mathbb F_{p^t}$.
Since $(\deg f_\lambda(x))^4=d^4<p^t$, $f_\lambda(x)$ is an exceptional polynomial over $\mathbb F_{p^t}$ from \cite[Theorem 8.4.19]{MP}.
Also, from \cite[Proposition 2.4]{BGQZ2017}, $f_\lambda(x)$ is indecomposable as exceptional polynomial over $\mathbb F_{p^t}$.
Hence, from \cite[Theorem 8.4.11]{MP}, $\deg f_\lambda(x)=2p^s+1$ is a prime not dividing $p^t-1$.
From the Niederreiter-Robinson criterion \cite[Lemma 1]{NR},
the polynomial $x^{\frac{(p^t)^{2p^s}-1}{p^t-1}+1}+bx$ is a PP of $\mathbb F_{(p^t)^{2p^s}}$; equivalently, the monomial $ b^{-1}x^{\frac{p^{2tp^s}-1}{p^t-1}+1} $ is a CPP of $\mathbb F_{p^{2tp^s}}$.
Thus, from \cite[Theorem 3.1]{BGQZ2017}, one of the following cases hold, where $\zeta$ is a primitive $(2p^s+1)$-th root of unity and $\alpha:=\zeta+\zeta^{-1}$, $\beta:=\zeta-\zeta^{-1}$:

\begin{itemize}

\item $p^t$ has order $2p^s$ modulo $2p^s+1$ and $b=\zeta-1$ up to multiplication by a non-zero element in $\mathbb F_{p^t}$. Since $b\in\mathbb F_{p^{2t}}\setminus\mathbb F_{p^t}$ and $2p^s+1$ is prime, we have $(2p^s+1)\mid(p^t+1)$. Hence $p^t\equiv-1\pmod{2p^s+1}$ has order $2\ne 2p^s$ modulo $2p^s+1$, a contradiction.

\item $p^t$ has order $2p^s$ modulo $2p^s+1$ and $b=e(\alpha-1)\sqrt{\beta^2(e^2-4a)}$ up to multiplication by a non-zero element in $\mathbb F_{p^t}$, for some $e,a\in\mathbb F_{p^t}$, $a\ne0$, such that $e^2-4a$ is a square in $\mathbb F_{p^t}$.
Since $b\in\mathbb F_{p^{2t}}\setminus\mathbb F_{p^t}$ we have $\alpha\in\mathbb F_{p^t}$, which implies $\zeta^{p^t-1}=1$ or $\zeta^{p^t+1}=1$. As $2p^s+1$ is prime and $b\notin\mathbb F_{p^t}$, this yields $(2p^s+1)\mid(p^t+1)$; hence, $p^t$ has order $2$ modulo $2p^s+1$, a contradiction to $s>0$.

\item $p^t$ has order $p^s$ modulo $2p^s+1$ and $b=e(\alpha-1)\sqrt{\beta^2(e^2-4a)}$ up to multiplication by a non-zero element in $\mathbb F_{p^t}$, for some $e,a\in\mathbb F_{p^t}$, $a\ne0$, such that $e^2-4a$ is $0$ or a non-square in $\mathbb F_{p^t}$.
From $b\in\mathbb F_{p^{2t}}\setminus\mathbb F_{p^t}$ we have $\alpha\in\mathbb F_{p^t}$ which implies $\zeta^{p^t-1}=1$ or $\zeta^{p^t+1}=1$; hence, $p^t$ has order $1$ or $2$ modulo $2p^s+1$, a contradiction to $s>0$.
\end{itemize}
Therefore, $f_\lambda(x)$ is not a PP of $\mathbb F_{p^t}$.

\item Case $r>1$.

Assume first $u\leq s$, so that we can take $u=0$ and $d=2p^s+r$. Suppose by contradiction that $f_\lambda(x)$ is a PP of $\mathbb F_{p^t}$. As $(\deg f_\lambda(x))^4<p^t$, $f_\lambda(x)$ is exceptional over $\mathbb F_{p^t}$, see \cite[Theorem 8.4.19]{MP}.
Note that $f_\lambda(x)$ has exactly three distinct zeros, one in $\mathbb F_{p^t}$ with multiplicity $r$ and two in $\mathbb F_{p^{2t}}\setminus\mathbb F_{p^t}$ with multiplicity $p^s$.

\begin{itemize}
\item
Suppose that $f(x)$ is indecomposable as exceptional polynomial over $\mathbb F_{p^t}$.
From \cite[Theorem 8.4.10]{MP}, one of the following cases holds.
\begin{itemize}
\item 

$2p^s+r=p^w$ for some $w\geq1$. In this case
$$f(x)=\left(x^{p^{w-s}}+x^{p^{w-s}-1}+\lambda x^{p^{w-s}-2}\right)^{p^s};$$
since $x\mapsto x^{p^s}$ is a permutation of $\mathbb F_{p^t}$, we can assume $s=0$. Then
\begin{small}
\begin{equation}\label{curve}
\frac{f(x)-f(y)}{x-y} = (x-y)^{p^w-1} + x^{p^w-2}+x^{p^w-3}y+\cdots+y^{p^w-2}+\lambda( x^{p^w-3}+x^{p^w-4}y+\cdots+y^{p^w-3} ).
\end{equation}
\end{small}
Let $\mathcal C$ be the plane curve of degree $p^w-1$ defined over $\mathbb F_{p^t}$ by affine equation $\frac{f(x)-f(y)}{x-y}=0$.
From Equation \eqref{curve}, $\mathcal C$ has a unique point at infinity $P_\infty$. Moreover, $\mathcal C$ intersects the line $x=y$ at the affine points $(0,0)$ and $(-2\lambda,-2\lambda)$ with multiplicity $p^w-3$ and $1$, respectively; hence, $P_\infty$ is a simple point for $\mathcal C$.
This implies that $\mathcal C$ is absolutely irreducible, a contradiction to the exceptionality of $f(x)$ (see \cite[Theorem 8.4.4]{MP}).

\item $2p^s+r=\frac{p^a(p^a-1)}{2}$, with $p\in\{2,3\}$ and $a>1$ odd; this is not possible, since $p\nmid r$.
\item $2p^s+r$ is coprime with $p$. From \cite[Theorem 8.4.11]{MP}, one of the following holds:
\begin{itemize}
\item $f_\lambda(x)$ is linear. This is not possible by the assumptions.
\item $f_\lambda(x)=x^{2p^s+r}$ where $2p^s+r$ is a prime not dividing $p^t-1$, up to composition with linear functions.
Then $f_\lambda(x)$ has either one or $n$ distinct roots in $\overline{\mathbb{F}}_{p^t}$, a contradiction.
\item $f_\lambda(x)=\ell_1\circ D_{2p^s+r}(\ell_2(x),a)$, where $2p^s+r$ is a prime not dividing $p^{2t}-1$, $D_{2p^s+r}(x,a)$ is a Dickson polynomial with $a\ne0$ of degree $2p^s+r$, and $\ell_1,\ell_2\in\mathbb F_{p^t}[x]$ are linear permutations.
If $(2p^s+r)\nmid(p^{2t}+1)$, then $D_{2p^s+r}(x,a)$ is a PP of $\mathbb F_{p^{2t}}$; see \cite[Theorem 8.4.11]{MP}. This is not possible, as $f_\lambda(x)$ has three distinct zeros in $\mathbb F_{p^{2t}}$.
Thus, $(2p^s+r)\mid(p^{2t}+1)$. Denote $\ell_1(x)=bx+c$. As $\ell_2$ permutes $\mathbb F_{p^{2t}}$, the number of zeros of $f_\lambda(x)$ in $\mathbb F_{p^{2t}}$ is equal to the number $Z$ of preimages of $-c/b$ under $D_{2p^s+r}(x,a)$; hence, $Z=3$.
On the other hand, from \cite[Theorems $3.26$ and $3.26^\prime$]{LMT1993} we have $Z\in\{1,2p^s+r,\frac{2p^s+r}{2},\frac{2p^s+r+1}{2}\}$.
Then $s=0$ and $r=3$, so that $f_\lambda(x)=x^5+x^4+\lambda x^3$ with $p\ne5$.
We have $D_5(x,a)=x^5-5ax^3+5a^2x$; by direct inspection, the polynomial $\ell_1\circ D_5(\ell_2(x),a)$ cannot have the form $x^5+x^4+\lambda x^3$ for any $\ell_1,\ell_2$.
\end{itemize}

\end{itemize}
\item
Now suppose that $f_\lambda(x)$ is a decomposable exceptional polynomial over $\mathbb F_{p^t}$, say $f_\lambda(x)=h(k(x))$ for some exceptional polynomials $h,k\in\mathbb F_{p^t}[x]$ with $\deg(h),\deg(k)>1$.
The roots of $f_\lambda(x)/x^r$ are conjugated under the Frobenius map $x\mapsto x^{p^t}$; hence, the polynomial
$$ \frac{f_\lambda(-x)-f_\lambda(0)}{(-x)^r}=\frac{h(k(-x))-h(k(0))}{(-x)^r}=\frac{(k(-x)-k(0))\prod_{i=1}^{\deg(h)-1}(k(-x)-\beta_i)}{(-x)^r} $$
is a power of a unique irreducible factor over $\mathbb F_{p^t}$.

Suppose that $k(-x)-k(0)$ has a monic absolutely irreducible factor different from $x$ and defined over $\mathbb F_{p^t}$.
Since the roots of $f_\lambda(x)/x^r$ are conjugated under $x\mapsto x^{p^t}$, we have $\beta_i=k(0)$ for all $i$. Hence, $\frac{f_\lambda(-x)-f_\lambda(0)}{(-x)^r}=\frac{(k(-x)-k(0))^{\deg(h)}}{(-x)^r}$. Also, $m\cdot\deg(h)=r$, where $x^m$ is the maximum power of $x$ which divides $k(-x)-k(0)$; in particular, $p\nmid\deg(h)$.
On the other hand, $f_\lambda(x)$ has just two distinct non-zero roots (the ones of $x^2+x+\mu$ where $\mu^{p^s}=\lambda$) with multiplicity $p^s$; hence, $\deg(h)\mid p^s$. This is a contradiction, either to $p\nmid\deg(h)$ or to $\deg(h)>1$.

Suppose that $k(-x)-k(0)=a x^m$, for some $a\in\mathbb F_{p^t}$ and $m>1$ with $\gcd(m,p^t-1)=1$.
If $p\mid m$, then $f_\lambda(x)$ is invariant under $x\mapsto \gamma x$ when $\gamma\in\mathbb F_p$; this is a contradiction to $\gcd(2p^s+r,p^s+r)=1$.
Then $p\nmid m$. Let $\bar x$ be a non-zero root of $f_\lambda(x)$; for any $\delta$ with $\delta^m=1$, $k(\delta\bar x)=k(\bar x)$ and $f_\lambda(\delta\bar x)=0$. Thus, the number of distinct non-zero roots of $f_\lambda(x)$ is a multiple of $m$; hence, $m=2$.
This implies $p=2$. Therefore $f_\lambda(x)=h(k(0)+ax^2)=h((\ell_0+\ell_1 x)^2)$ with $\ell_0,\ell_1\in\mathbb F_{p^t}$, so that the polynomial $h(x^2)$ is also exceptional of degree $\deg(f_\lambda)$. Since $\deg(h)$ is odd, this is not possible.
\end{itemize}

We have shown that $f_\lambda(x)$ is not a PP of $\mathbb F_{p^t}$ under the assumption $u\leq s$.
If $u>s$, then we can take $s=0$ so that $d=r+2$ and $f_\lambda(x)=x^r(x^2+x+\lambda)$. The same arguments as in the case $u\leq s$ still apply and show that $f_\lambda(x)$ is not a PP of $\mathbb F_{p^t}$.
\end{itemize}

\begin{remark}
Theorem {\rm \ref{result}} yields the characterization also of permutation trinomials of $\mathbb F_{p^t}$ of type $g_{\alpha,\beta}(x)=x^{2p^s+r}+\alpha x^{p^s+r}+\beta x^r$, under the assumptions $\alpha\ne0$ and $d^4<p^t$ (with $d$ defined as in Theorem {\rm \ref{result}}).

In fact, let $\gamma\in\mathbb F_{p^t}$ satisfy $\gamma^{p^s}=\alpha$. Then $g_{\alpha,\beta}(x)$ is a PP of $\mathbb F_{p^t}$ if and only if $\frac{1}{\gamma^{2p^s+r}}g(\gamma x)=f_{\beta/\alpha^2}(x)$ is a PP of $\mathbb F_{p^t}$. Thus, $g_{\alpha,\beta}(x)$ is a PP of $\mathbb F_{p^t}$ exactly in the following cases:
\begin{itemize}
\item $p=2$, $t$ is odd, and $g_{\alpha,\beta}(x)=x^3+\alpha x^2+\alpha^2 x$ or $g_{\alpha,\beta}(x)=x^5+\alpha x^3+\alpha^2 x$ for some $\alpha\in\mathbb F_{p^t}$;
\item $p\equiv2\pmod3$, $t$ is odd, and $g_{\alpha,\beta}(x)=x^3+\alpha x^2+\frac{\alpha^2}{3}x$ for some $\alpha\in\mathbb F_{p^t}$.
\end{itemize}
\end{remark}


\begin{flushleft}
Daniele Bartoli\\
Department of Mathematics and Computer Science,\\
University of Perugia,\\
e-mail: {\sf daniele.bartoli@unipg.it}
\end{flushleft}

\begin{flushleft}
Giovanni Zini\\
Department of Mathematics and Computer Science ``Ulisse Dini'',\\
University of Florence,\\
e-mail: {\sf gzini@math.unifi.it}
\end{flushleft}

\end{document}